\documentclass[12pt,a4paper,reqno]{amsart}
\usepackage[margin=2cm]{geometry}
\usepackage{graphics, amsmath,amssymb, enumitem, comment, url, mathtools}
\usepackage{graphicx}
\usepackage{epsfig}

\newif\ifcolorcomments
\newcommand{\allowcomments}[4]{
\newcommand{#1}[1]{\ifdraft{\ifcolorcomments{\textcolor{#4}{##1 --#3}}\else{\textsl{ ##1 \ --#3}}\fi}\else{}\fi}
}

\allowcomments{\commumtaz}{MH}{Mumtaz}{green}
\allowcomments{\comjohannes}{JS}{Johannes}{blue}
\allowcomments{\comdavid}{DS}{DS}{magenta}
\allowcomments{\comBW}{BW}{BW}{red}

\colorcommentstrue
\usepackage{amssymb}
\usepackage{amsmath,amsthm}

\usepackage{colortbl}

\newtheorem{theorem}{Theorem}[section]

\newtheorem{proposition}[theorem]{Proposition}

\theoremstyle{definition}

\newtheorem{remark}[theorem]{Remark}

\newcommand{\FF}{\mathcal F}

\newcommand{\HH}{\mathcal H}

\newcommand{\MM}{\mathcal M}

\newcommand{\N}{\mathbb N}

\newcommand{\R}{\mathbb R}

\newcommand{\Z}{\mathbb Z}

\newcommand{\mbf}{\mathbf}

\newcommand{\pp}{\mbf p}

\newcommand{\xx}{\mbf x}

\renewcommand{\text}{\textup}

\newcommand{\NPC}[1]{\ignorespaces}

\newcommand{\bftheta}{{\boldsymbol\theta}}
\newcommand{\bftau}{{\boldsymbol\tau}}

\newif\ifdraft\drafttrue
\def\N{\mathbb N}
\def\Z{\mathbb Z}

\def\R{\mathbb R}

\IfFileExists{marvosym.sty}{
\RequirePackage{marvosym} 
}{

}

\IfFileExists{wasysym.sty}{
\RequirePackage{wasysym}\renewcommand{\emptyset}{{\diameter}}
}{

}

\usepackage{setspace}
\newcommand*{\myDots}{\ifmmode\mathellipsis\else.\kern-0.07em.\kern-0.07em.\fi}

\begin{document}
\title{An exponentially shrinking problem}

\author{Mumtaz ~Hussain}
\author{Junjie Shi}
\address{Mumtaz Hussain, Department of mathematical and physical sciences, La Trobe University,  Bendigo 3552, Australia. }
\email{m.hussain@latrobe.edu.au}
\address{Junjie Shi, School of Mathematics and Statistics, Huazhong University of Science and Technology, Wuhan 430074, China}
\email{junjie\_Shi@hust.edu.cn}

\begin{abstract} The Jarn\'ik-Besicovitch theorem is a fundamental result in metric number theory which concerns the Hausdorff dimension for some limsup sets. We investigate a related problem of estimating the Hausdorff dimension of a liminf set.  Let $\tau>0$ and $\{q_j\}_{j\ge 1}$ be a sequence of integers.  We calculate the Hausdorff dimension of the set $$\Lambda^\bftheta_d(\tau)=\left\{\xx\in[0, 1)^d: \|q_jx_i-\theta_i\|<q_j^{-\tau} \ \text{for all } j\geq 1, i=1,2,\cdots,d\right\},$$
where $\left\|\cdot\right\|$ denotes the distance to the nearest integer and $\bftheta\in [0, 1)^d$ is fixed.  We also give some heuristics for the Hausdorff dimension of the corresponding multiplicative set
$$\MM_d^\bftheta(\tau)=\left\{\xx\in[0, 1)^d:\prod_{i=1}^d \|q_jx_i-\theta_i\|<q_j^{-\tau} \ \text{for all } j\geq 1\right\}.$$

\end{abstract}
\maketitle

\section{Introduction}\label{sec1} How well a real number can be approximated by rational numbers is a fundamental problem in metric number theory. Dirichlet's theorem (1842) states that for any $\xx:=(x_1, \ldots, x_d)\in [0, 1)^d$ and any $Q\in\N$, there exixts $1\leq q\leq Q$ such that
\begin{equation*}\label{DTeqn}
\|qx_i\|\leq Q^{-\frac{1}{d}}, \quad 1\leq i\leq d,
\end{equation*}
where $\left\|\cdot\right\|$ denotes the distance to the nearest integer. A consequence of Dirichlet's theorem is that for $\xx\in[0, 1)^d$, there exists infinitely many $q\in\N$ such that
\begin{equation}\label{KTeqn}
\|q\xx\|:=\max\left\{\|qx_1\|, \ldots, \|qx_d\|\right\}\leq q^{-\frac{1}{d}}.
\end{equation}
Naturally, one can ask what happens when the right hand side of \eqref{KTeqn} is replaced by some general function $\psi:\N\mapsto \R_+$ such that $\psi(q)\to 0$ as $q\to\infty$. This setup, called the set of $\psi$-approximable numbers, is defined as
$$W_d(\psi)=\left\{\xx\in[0, 1)^d: \|q\xx\|<\psi(q) \ \text{for infinitely many } q\in\N\right\}.$$
Results on the metrical size of $W_d(\psi)$ are well-known. These include fundamental results such as Khintchine's  theorem (1924) and Jarn\'ik's  theorem (1932) on the Lebesgue measure and Hausdorff measure of $W_d(\psi)$ respectively.  For recent advances in this area, we refer the reader to the recent breakthrough result of Maynard-Koukouloupolous \cite{KMDS} for the measure of $W_1(\psi)$ without the monotonicity assumption on $\psi$.

It should be clear that any point $\xx$ in the limsup set $W_d(\psi)$ can be well approximated by infinitely many rational points, however there is no further information about the rational numbers approximating $\xx$, the rational numbers are rather generic.  So a natural question is ``What is the size of the set $W_d(\psi)$ if {\em infinitely many} is replaced by {\em for all} within the set $W_d(\psi)$\footnote{The authors thanks Johannes Schleischitz for raising this question. In particular, he identified a connection between this question and Theorems 2.3 and 2.4 in his recent paper \cite{Schleischitzuniform} regarding the exact approximation and Dirichlet spectrum. }?''

To the best of our knowledge, any explicit metrical result is not known for this set.
In this paper, we consider this problem for the  inhomogeneous settings that encompasses the homogeneous settings discussed above.

Let $\{q_j\}_{j\geq1}$ be a sequence of positive integers.
Define the set
$$\Lambda^\bftheta_d(\tau)=\left\{\xx\in[0, 1)^d: \|q_jx_i-\theta_i\|<q_j^{-\tau} \ \text{for all } j\geq 1, i=1,2,\cdots,d\right\}.$$
%
%
%
%
%
%
%
%
%
%
For the definition of Hausdorff measure and dimension we refer the reader to \cite{Falconer_book2003}. Throughout, we use the notation  $B(c, r)$ to denote a ball centred at $c$ with radius $r$. For real quantities $A, B$ that
depend on parameters, we write $A\ll B$ if $A\leq c B$ for a constant $c > 0$ that is independent of those
parameters. We write $A\asymp B$ if $A\ll B\ll A$. 

%

%

\begin{theorem}\label{thm2} Assume $h> \tau+1$. Then
$$\dim_\HH \Lambda^\bftheta_d(\tau)=\frac{d(1-\tau\alpha)}{\tau+1},$$ 
where
\begin{equation}\label{liminf}h=\liminf_{j\to\infty}\frac{\log q_{j+1}}{\log q_j}, \ \ \ \alpha=\liminf_{j\to \infty}\frac{\log q_1+\cdots+\log q_{j-1}}{\log q_j}.
\end{equation}
\end{theorem}

\begin{remark}\label{r1} When $h\leq \tau+1$, the structure of the set $\Lambda^\bftheta_d(\tau)$ is quite complex and, hence, poses a challenging question of the Hausdorff dimension of the $\Lambda^\bftheta_d(\tau)$ for this case. We present two examples to highlight the complexity for the one dimensional homogeneous case only. Fix a real number $\tau>0$ and let $\theta_j=0$ for all $j\ge 1$.

\begin{itemize}
\item[(1)]Assume that the integer sequence $\{q_j\}_{j\ge 1}$ satisfies$$
\frac{1}{8}q_j^{\tau+1}\le q_{j+1}\le  \frac{1}{4}q_j^{1+\tau}, \ \ {\text{for all}}\ j\ge 1.
$$ Then $h=1+\tau$. We claim that $\Lambda^\bftheta_d(\tau)$
contains at most finitely many points.

More precisely, for each $j\ge 1$, let $$E_j=\{x\in[0, 1): \|q_jx\|<q_j^{-\tau} \}.$$ Then $E_{j}$ consists of $q_{j}$ intervals of maximum length $2q_j^{-1-\tau}$. At the same time, the intervals contained in $E_{j+1}$ are separated by a gap at least $q_{j+1}^{-1}-2q_{j+1}^{-1-\tau}$.
Note that if $q_1$ is large enough then for all $j\ge 1$$$
\frac{2}{q_j^{1+\tau}}\leq \frac{1}{q_{j+1}}-\frac{2}{q_{j+1}^{1+\tau}}.
$$ This implies that each interval inside $E_j$ can intersect at most one interval in $E_{j+1}$. Thus $E_j\cap E_{j+1}$ will be contained in at most $q_{j}$ intervals of length $2q_{j+1}^{-1-\tau}$. Continuing in this way, we see that the set $$
E_1\cap E_2\cap\cdots \cap E_j
$$ will be contained in at most $q_1$ intervals of length $2q_j^{-1-\tau}$. Taking the limit as $j$ tends to infinity shows the claim.

\item[ (2)] Let $\eta>1+\tau$. Define a sequence of integers $\{q_j\}_{j\ge 1}$ in the following way: $$
 \left\{
   \begin{array}{ll}
    \liminf\limits_{j\to\infty}\frac{\log q_{j+1}}{\log q_j}=1+\tau, \ q_{j+1}\le \frac{1}{4}q_{j}^{1+\tau}, \ q_j|q_{j+1}, & \hbox{when $j$ is even;} \\
    \liminf\limits_{j\to\infty}\frac{\log q_{j+1}}{\log q_j}=\eta, & \hbox{when $j$ is odd.}
   \end{array}
 \right.
 $$


With the same argument as in the first example, we note that for all $j\ge 1$, 
\begin{align*}
E_{2j}\cap E_{2j+1}&=\left\{x\in [0,1): \|q_{2j}x\|<q_{2j+1}^{-\tau}\right\}\\ &\approx\left\{x\in [0,1): \|q_{2j}x\|<q_{2j}^{-\tau(2+\tau)}\right\}.
\end{align*}
By taking $\bar{q}_j=q_{2j}$, we have$$
\Lambda^\bftheta_d(\tau)=\Big\{x\in [0,1): \|\bar{q}_{j}x\|<(\bar{q}_{j})^{-\tau(2+\tau)}, \ {\text{for all}}\ j\ge 1\Big\}.
$$ This falls into the case that $h>1+\tau$ for the sequence $\{\bar{q}_j\}_{j\ge 1}$ by noting that $h=\eta(1+\tau)$, and our new $\hat{\tau}\mapsto \tau(2+\tau)$. So $h=\eta(1+\tau)>(1+\tau)^2=\hat{\tau}+1$. Thus,  by Theorem \ref{thm2} we note that the set  $\Lambda^\bftheta_d(\tau)$ 
is of positive Hausdorff dimension.
\end{itemize}

%
%
%
\end{remark}
%

%


%
%

\section{Proof of Theorem \ref{thm2}}
We split the proof into two parts: the upper and lower bounds. For all $j\ge 1$, let $$E_j=\{\xx\in[0, 1)^d: \|q_jx_i-\theta_i\|<q_j^{-\tau} \ i=1,2,\cdots,d\}. $$
Then, $$\Lambda^\bftheta_d(\tau)=\bigcap_{j\ge 1}E_j.$$
By viewing $[0,1]^d$ as the unit torus, it is clear that $E_j$ is a collection of squares: 
\begin{equation*}\label{f1}
E_j= \bigcup_{0\le p_{i,j}<q_j, 1\le i\le d}\ \prod_{i=1}^dB\left(\frac{p_{i,j}+\theta_i}{q_{j}}, \frac{1}{q_j^{1+\tau}}\right).
\end{equation*} The following simple fact will be used in the upper and lower bound estimates below.

{\em Fact: Let $(a,b)$ be an interval in $\mathbb{R}$. For any $\theta\in [0,1)$ and $q\in \N$, the interval $(a,b)$ can contain $$
{\text{at most}}\ (b-a)\cdot q+2, \ \ \ {\text{and at least}}\ (b-a)q-2
$$ shifted rationals $(p+\theta)/q$.}

\subsection{The upper bound} It is trivial that $E_1\cap \cdots \cap E_j$ is a cover for $\Lambda^\bftheta_d(\tau)$. Next we search for a suitable cover for $E_1\cap \cdots \cap E_j$.

To elaborate the main strategy, we first consider $E_1\cap E_2$. Note that $E_1$ consists of $q_1^d$ squares of side length $2q_1^{-1-\tau}$, and the squares in $E_2$ are separated by a gap of at least $q_2^{-1}-q_2^{-1-\tau}\ge 1/(2q_2)$ in each coordinate axis. Thus each square in $E_1$ can intersect at most $$
\Big(2q_1^{-1-\tau}\cdot 2q_2+2\Big)^d
$$ squares in $E_2$. Hence, {by the {\em Fact} stated at the beginning,} $E_1\cap E_2$ can be covered by $$
q_1^d\cdot \Big(2q_1^{-1-\tau}\cdot 2q_2+2\Big)^d
$$ squares of side length $2q_2^{-1-\tau}$.

Thus by an iterative process, $E_1\cap \cdots \cap E_j$ can be covered by at most
$$q_1^d\cdot \Big({ 2}q_1^{-1-\tau}\cdot 2q_2+2\Big)^d\cdots \Big({ 2}q_{j-1}^{-1-\tau}\cdot 2q_j+2\Big)^d$$
squares of side length $2q_{j}^{-1-\tau}$.

Recall the definition of $h$ and the assumption $h> 1+\tau$. For any $\epsilon>0$, there exists $k_o\in \N$ such that for all $k\ge k_o$, $$
4q_k^{-1-\tau}\cdot q_{k+1}+2\le q_k^{-1-\tau}\cdot q_{k}^{1+\epsilon}.
$$
As a result, for some constant $c=c(\epsilon)$, one has \begin{align*}
  q_1^d\cdot \Big(q_1^{-1-\tau}\cdot 2q_2+2\Big)^d\cdots \Big(q_{j-1}^{-1-\tau}\cdot 2q_j+2\Big)^d\le \Big[c(\epsilon)\cdot q_1^{-\tau+\epsilon}\cdots q_{j-1}^{-\tau+\epsilon}\cdot q_j^{1+\epsilon}\Big]^d.
\end{align*}
Thus,
$$\HH^s(\Lambda^\bftheta_d(\tau))\leq \liminf_{j\to\infty}\Big[c(\epsilon)\cdot q_1^{-\tau+\epsilon}\cdots q_{j-1}^{-\tau+\epsilon}\cdot q_j^{1+\epsilon}\Big]^d\cdot\left(q_{j}^{-s(\tau+1)} \right).$$ Therefore $\HH^s(\Lambda^\bftheta_d(\tau))=0,$ whenever
\begin{align*}s>&\frac{d}{\tau+1}\liminf_{j\to \infty}\frac{(-\tau+\epsilon)(\log q_1+\cdots+\log q_{j-1})+(1+\epsilon)\log q_j}{\log q_j}\\
=&\frac{d}{\tau+1}\Big(1+\epsilon-(\tau-\epsilon)\alpha\Big).
\end{align*} Hence, from the definition of Hausdorff dimension and the arbitrariness of $\epsilon>0$, it follows that$$\dim_\HH \Lambda^\bftheta_d(\tau)\leq \frac{d(1-\tau\alpha)}{\tau+1}.$$

\subsection{The lower bound} 

The main tool in establishing the lower bound for the Hausdorff dimension of $ \Lambda^\bftheta_d(\tau)$ is the following well-known mass distribution principle \cite{Falconer_book2003}.
\begin{proposition}[Mass Distribution Principle \cite{Falconer_book2003}]\label{p1}
Let $\mu$ be a probability measure supported on a measurable set $F$. Suppose there are positive constants $c$ and $r_0$ such that
$$\mu(B(x,r))\le c r^s$$
for any ball $B(x,r)$ with radius $r\le r_0$ and center $x\in F$. Then $\dim_\HH F\ge s$.
\end{proposition}
To use this principle, we construct a suitable Cantor subset of $\Lambda^\bftheta_d(\tau)$ satisfying conditions of the mass distribution principle.

\begin{itemize}
\item Level one, $F_1$, of the Cantor subset is nothing but $E_1$, that is $$
F_1= \bigcup_{0\le p_{i,1}<q_1, 1\le i\le d}\ \prod_{i=1}^dB\left(\frac{p_{i,1}+\theta_i}{q_{1}}, \frac{1}{q_1^{1+\tau}}\right).
$$ We also write
$$\FF_1=\left\{\prod_{i=1}^d B\left(\frac{p_{i,1}+\theta_i}{q_1}, \tfrac1{q_1^{\tau+1}}\right): 0\leq p_{i,1}< q_1, \ 1\le i\le d\right\},$$ for the collection of squares constituting $E_1$.

\item For level two of the Cantor subset we first count the number of squares in $E_2$ contained in each square of $\FF_1.$ Let $B_1$ be a square in $\FF_1$ and write it as $$
B_1=\prod_{i=1}^dB\left(\frac{p_{i,1}+\theta_i}{q_{1}}, \frac{1}{q_1^{1+\tau}}\right):=\prod_{i=1}^dB_{i,1}.
$$ {by the {\em Fact} stated at the beginning,} for any $i\in\{1,\cdots,d\},$
    \begin{align*}\# P_{i,2}:&=\#\left\{0\le p_{i,2}<q_2: B\left(\frac{p_{i,2}+\theta_i}{q_2}, \frac1{q_2^{\tau+1}}\right)\cap B_{i,1}\neq\emptyset\right\}\\ &\geq\frac{2q_2}{q_1^{\tau+1}}-2\geq\frac{q_2}{q_1^{\tau+1}}.
    \end{align*}
    Thus, there are at least $[{q_2}\cdot {q_1^{-(\tau+1)}}]^d$ squares contained within each square of $\FF_1.$ Choose a sub-collection of these squares with cardinality $[{q_2}\cdot {q_1^{-(\tau+1)}}]^d$ and denote this sub-collection by $\FF_2(B_1)$. Then define $$
\FF_2=\bigcup_{B_1\in \FF_1}\FF_2(B_1), \ F_2=\bigcup_{B_2\in \FF_2}B_2,
$$as the second level of the Cantor subset.
Note that every pair of squares in $E_2$ is separated by a gap of at least $\tfrac{1}{2q_2},$ so are those in $\FF_2$.
This is because for any integer vectors
$(p_{1,2},\cdots, p_{d,2})\ne (p'_{1,2},\cdots, p'_{d,2})$, at least for one index $i$,
    \begin{align*}
      \left|\prod_{i=1}^d B\left(\tfrac{p_{i,2}+\theta_i}{q_2},\tfrac{1}{q_2^{\tau+1}}\right),\prod_{i=1}^d B\left(\tfrac{p'_{i, 2}+\theta_i}{q_2},\tfrac{1}{q_2^{\tau+1}}\right)\right|
        &=\max_{1\leq i\leq d} \left| B\left(\tfrac{p_{i, 2}+\theta_i}{q_2},\tfrac{1}{q_2^{\tau+1}}\right),B\left(\tfrac{p'_{i, 2}+\theta_i}{q_2},\tfrac{1}{q_2^{\tau+1}}\right)\right| \\
      &\geq\max_{1\leq i\leq d}\tfrac{p_{i, 2}-p'_{i, 2}}{q_2}-\tfrac{2}{q_2^{\tau+1}}\\ &\geq\tfrac{1}{q_2}-\tfrac{2}{q_2^{\tau+1}}\\ &\geq\tfrac{1}{2q_2}.
    \end{align*}
\end{itemize} 
Continuing in this way, we have a sequence $\{\FF_k\}_{k\ge 1}$ consisting of squares such that $$
\FF_k=\bigcup_{B_{k-1}\in \FF_{k-1}}\FF_k(B_{k-1}).
$$ For each $B_{k-1}\in \FF_{k-1}$, the family $\FF_k(B_{k-1})$ contains $[{q_k}\cdot {q_{k-1}^{-(\tau+1)}}]^d$ squares of side length $q_{k}^{-(1+\tau)}$ and separated by at least $\tfrac{1}{2q_k}$ gap between them. Then define $$
\FF_k=\bigcup_{B_k\in \FF_k}B_k, \ {\text{and}}\ \FF_{\infty}=\bigcap_{k\ge 1}F_k,
$$as the $k$th level of the Cantor subset and the Cantor subset of $\Lambda^\bftheta_d(\tau)$ respectively.

Next, we distribute the mass uniformly on the Cantor subset.  For any $B_k\in\FF_k,$ let $B_{k-1}$ be the square in $\FF_{k-1}$ such that $B_k\in \FF_k(B_{k-1})$ (when $k=1$, define $B_0=[0,1]^d$). Then define
\begin{align*}\mu(B_k)&=\mu(B_{k-1})\frac1{\#\FF_k(B_{k-1})}\\ &=\mu(B_{k-1})\left[\frac{q_k}{q_{k-1}^{\tau+1}}\right]^{-d}
\\ &=\left(\prod_{i=1}^k\left[\frac{q_{i-1}^{\tau+1}}{q_i}\right]\right)^d\\ &\asymp{q_k^{-d}\cdot (q_1\cdots q_{k-1})^{d\tau}}.
\end{align*}

Having defined the measure on the Cantor subset, we now calculate the H\"{o}lder exponent of the measure $\mu$, which determines how the measure of a ball is related to its diameter. At first, we consider the measure of balls in $\FF_k$: \begin{align*}
  \liminf_{k\to\infty}\frac{\log \mu(B_k)}{|B_k|}&=\liminf_{k\to\infty}\frac{d\log q_k-d\tau(\log q_1+\cdots+\log q_{k-1})}{(\tau+1)\log q_k}\\ &=\frac{d(1-\alpha)}{\tau+1},
\end{align*}
where $|B_k|$ denotes the radius of the ball $B_k$. Thus for any $s<\frac{d(1-\alpha)}{\tau+1}$, there exists $k_o\in\N$ such that for all $k\ge k_o$
$$
\mu(B_k)\le |B_k|^{s},\ \ {\text{for all}}\ B_k\in\FF_k.
$$

For any $x\in \FF_\infty$ and $r<1/(4q_{k_o})$, we consider the measure of the ball $B(x,r)$. Assuming the ball $B(x,r)$ can intersect only one element in $\FF_{k-1},$ say $B_{k-1},$ and at least two elements in $\FF_k.$ Hence all the elements in $\FF_k$ which can intersect $B(x,r)$ are contained in $\FF_k(B_{k-1}).$ Since any two elements in  $\FF_k(B_{k-1})$ are $\tfrac1{2q_k}$-separated, we have $r\geq\tfrac1{4q_k}.$ Thus $k\ge k_o$.

Without loss of generality, we can assume that $r\leq|B_{k-1}|,$ otherwise, we have
$$\mu (B(x,r))\leq\mu (B_{k-1})\leq|B_{k-1}|^s\leq r^s.$$
We count how many elements in $\FF_k(B_{k-1})$ which can intersect $B(x,r).$ Note that all elements in $\FF_k(B_{k-1})$ are
$\tfrac1{2q_k}$-separated, which mean that their $\tfrac1{4q_k}$-thickenings are still disjoint. Moreover, since $r\geq\tfrac1{4q_k},$ all the $\tfrac1{2q_k}$-thickenings of these elements in $\FF_k(B_{k-1})$ are still contained in $B(x,8r).$ Thus, a volume argument gives $$\#\{B_k\in\FF_k(B_{k-1}): B_k\cap B(x,r)\neq\emptyset\}\leq 16^d\cdot r^d\cdot q_k^d.$$

As a consequence,
\begin{align*}
   \mu(B(x,r)) & \leq\min\left\{\mu(B_{k-1}),\ \ 16^d\cdot r^dq_k^d\left(\tfrac{q_k}{q_{k-1}^{\tau+1}}\right)^{-d}\mu(B_{k-1})\right\} \\
    &=\min\left\{\mu(B_{k-1}),\ 16^d\cdot r^dq_{k-1}^{d(\tau+1)}\mu(B_{k-1})\right\}\\
    &\leq 16^d\cdot \left({q_{k-1}^{-(\tau+1)}}\right)^s\min\left\{1,r^dq_{k-1}^{d(\tau+1)}\right\}\\
    &\leq16^d\cdot \left({q_{k-1}^{-(\tau+1)}}\right)^s\cdot 1^{(1-s/d)}\left(r^dq_{k-1}^{d(\tau+1)}\right)^{s/d}\\
    &=16^d\cdot r^s.
 \end{align*}
 The Mass Distribution Principle implies that
 $$\dim_\HH \Lambda^\bftheta_d(\tau)\geq \frac{d(1-\alpha)}{\tau+1}.$$
 This completes the proof.

\section{Concluding remarks and open problems}
Finally, we list a couple of open problems which are worth investigating. The first problem may be termed as the weighted inhomogeneous exponentially shrinking problem.  Fix $\bftheta:=(\theta_1, \ldots, \theta_d)\in [0,1)^d$ and $\bftau=(\tau_1, \ldots, \tau_d)\in\R_+^d$. Let the integer sequence $\{q_j\}_{j\geq1}$ be as in \eqref{liminf}. Then, what is the Hausdorff dimension of the set
$$\Lambda_d^{\bftheta}(\bftau)=\left\{x\in[0, 1)^d: \|q_jx_i-\theta_i\|<q_j^{-\tau_i} \ \text{for all } j\geq 1, i=1,2,\cdots,d\right\}?$$
We expect that the ideas presented in this paper combined with techniques in \cite{WWX} may yield a solution to this problem. A more general problem may be to estimate the Hausdorff measure of
$$\Lambda_d^{\bftheta}(\psi)=\left\{\xx\in[0, 1)^d: \|q_jx_i-\theta_i\|<\psi(q_j) \ \text{for all } j\geq 1, i=1,2,\cdots,d\right\}$$
for some function $\psi:\N\to \R_+$.

In another direction, the same problem may be considered within the multiplicative settings.   For any $\tau>0$ and fixed $\bftheta\in[0,1)^d$, consider the set
%
%
$$\MM_d^\bftheta(\tau)=\left\{\xx\in[0, 1)^d:\prod_{i=1}^d \|q_jx_i-\theta_i\|<q_j^{-\tau} \ \text{for all } j\geq 1\right\}.$$

We note that, since $\|q_jx_1-\theta_1\|<1$, we trivially have that $$ \MM_1^{\theta_1}(\tau)\times [0,1)^{d-1}\subseteq \MM_d^{\bf{\theta}}(\tau).$$ Hence giving the lower bound of the Hausdorff dimension to be $$d-1+ \dim_\HH\MM_1^\theta(\tau)=d-1+ \dim_\HH\Lambda_1^\theta(\tau)=d-1+ \tfrac{1-\tau\alpha}{\tau+1},$$
where $$\alpha=\liminf_{j\to \infty}\frac{\log q_1+\cdots+\log q_{j-1}}{\log q_j}, \quad \text{and}\quad \tau+1<h=\liminf_{j\to\infty}\frac{\log q_{j}}{\log q_{j-1}}.$$
So the main problem is to estimate the upper bound for $\dim_\HH \MM_d^{\bf{\theta}}(\tau).$ We prove the following theorem.
%

\begin{theorem}\label{thm3} Let $h> \tau+1$, then
$$d-1+ \tfrac{1-\tau\alpha}{\tau+1} \leq \MM_d^\bftheta(\tau)\leq d-1+\frac{1}{\tau+1}.$$
\end{theorem}

\subsection{Proof of Theorem \ref{thm3}}For the upper bound, we proceed as in \cite{HussainSimmons}.  Let $h>0$ and define the sequence of positive integers  as in \eqref{liminf}. Then
\begin{align*}
\MM_d^\bftheta(\tau)&=\left\{\xx\in[0, 1)^d:\prod_{i=1}^d \|q_jx_i-\theta_i\|<q_j^{-\tau} \ \text{for all } j\geq 1\right\}\\
&\subseteq \bigcap_{q\in\N}\bigcup_{p\in \mathcal Z_q} \tfrac1{q_j}(\pp+\bftheta+M(q_j^{-\tau})),
\end{align*}
where
\begin{align*}
\mathcal Z_{q_j}&=\{\pp\in \Z^d: -1\leq p_i+\theta_i\leq q_j+1 \quad \forall i=1,\ldots, d\}\\
M(\gamma)&=\{\xx\in[0, 1)^d:\prod_i|x_i|\leq \gamma\}.
\end{align*}
For each $q_j$, let $\{B_{{q_j}, i}: i=1, \ldots, N_{q_j}\}$ be the collection of $d$-dimensional hypercubes covering $M(q_j^{-\tau})$ as in \cite[Lemma 2.2]{HussainSimmons}.  Then the collection
$$\left\{\tfrac1{q_j}(\pp+\bftheta+B_{q_j, i}): q_j\in\N, \pp\in\mathcal Z_{q_j}\right\}$$
is a fine cover for $\MM_d^\bftheta(\tau)$.
Therefore, the $s$-dimensional cost of the cover can be estimated as
\begin{align*}
\HH^s(\MM_d(\tau))&=\liminf_{j\to\infty}\sum_{p\in\mathcal Z_{q_j}}\sum_{i=1}^{N_q}\left|\tfrac1{q_j}(\pp+\bftheta+B_{q_j, i})\right|^s
\\ &\ll \liminf_{j\to\infty} q_j^{d} q_j^{-s} q_j^{-\tau(s-d+1)}\\
&\ll \liminf_{j\to\infty}  q_j^{d-s-\tau(s-d+1)}.
\end{align*}
Hence, when $s>d-1+\tfrac{1}{\tau+1}$, $\HH^s(\MM_d(\tau))\to 0$. Thus
$$\dim_\HH \MM_d(\tau)\leq d-1+\tfrac{1}{\tau+1}.$$
Finally, note that the dimension bounds are not  sharp. A careful covering argument may lead to the dimension to be
$$\dim_\HH\MM_d^\bftheta(\tau)=d-1+ \tfrac{1-\tau\alpha}{\tau+1}. $$

Finally, it is worth noting that for $h$ large enough $\alpha\approx \frac1h$. So for $\tau$ large enough ($\tau+1\approx \tau$) with $\tau+1=h-\epsilon$, both the dimension bounds are close enough.

\medskip

\noindent{\bf Acknowledgements.} This research is supported by the Australian Research Council Discovery Project (200100994). The author thanks Dr Johannes Schleischitz for proposing this problem. We also thank Professor Yann Bugeaud, Professor Baowei Wang, and Dr Ben Ward for useful discussions. We thank an anonymous referee and the handling editor for several comments and suggestions for improvement.

%
%

\providecommand{\bysame}{\leavevmode\hbox to3em{\hrulefill}\thinspace}
\providecommand{\MR}{\relax\ifhmode\unskip\space\fi MR }
\providecommand{\MRhref}[2]{%
  \href{http://www.ams.org/mathscinet-getitem?mr=#1}{#2}
}
\providecommand{\href}[2]{#2}

\end{document}